\documentclass[12pt,reqno]{amsart}
\usepackage{amscd}
\input psfig.sty

\newcommand{\p}{\partial}
\newcommand{\ep}{\varepsilon}
\newcommand{\sign}{\operatorname{sign}}
\newcommand{\rank}{\operatorname{rank}}
\newcommand{\coker}{\operatorname{coker}}
\newcommand{\grad}{\operatorname{grad}}
\newcommand{\Hom}{\operatorname{Hom}}
\newcommand{\hol}{\operatorname{hol}}
\newcommand{\ad}{\operatorname{ad}}
\newcommand{\SU}{\operatorname{SU}}
\newcommand{\SO}{\operatorname{SO}}
\newcommand{\Fix}{\operatorname{Fix}}
\newcommand{\ind}{\operatorname{index}}
\newcommand{\su}{\mathfrak{su}}
\newcommand{\tr}{\operatorname{tr}}
\renewcommand{\Re}{\operatorname{Re}}
\renewcommand{\fbox}{}

\newtheorem{theorem}{Theorem}
\newtheorem{lemma}[theorem]{Lemma}
\newtheorem{corollary}[theorem]{Corollary}
\newtheorem{proposition}[theorem]{Proposition}

\title[Brieskorn homology spheres]
      {Floer homology of Brieskorn homology spheres : \\ 
       solution to Atiyah's problem} 
\author{Nikolai Saveliev}
\subjclass{57M25, 57R80}
\keywords{Floer homology, Seifert manifolds, Montesinos knots, Casson
invariant, Jones polynomial, homology cobordism, knot signature.} 
\thanks{Research at MSRI is supported in part by NSF grant DMS-9022140.}
\address{\hskip-\parindent Department of Mathematics, University of Michigan,
Ann Arbor, MI 48109} 
\email{saveliev@math.lsa.umich.edu}

\begin{document}
\begin{abstract} 
In this paper we answer the question posed by M.~Atiyah, see \cite{NW}, and 
give an explicit formula for Floer homology of Brieskorn homology spheres in 
terms of their branching sets over the 3--sphere. We further show how Floer 
homology is related to other invariants of knots and 3--manifolds, among which 
are the $\bar\mu$--invariant of W.~Neumann and L.~Siebenmann and the Jones 
polynomial. Essential progress is made in proving the homology cobordism 
invariance of our own $\nu$--invariant, see \cite{S}.
\end{abstract}

\maketitle
Let $p,q$, and $r$ be pairwise coprime positive integers. The Brieskorn
homology 3--sphere $\Sigma (p,q,r)$ is the link of the singularity of
$f^{-1}(0)$, where $f:{\mathbb C}^3 \to {\mathbb C}$ is a map of the form 
$f(x,y,z) = x^p+y^q+z^r$. The complex conjugation in ${\mathbb C}^3$ acts on 
$\Sigma (p,q,r)$ turning it into a double branched cover of $S^3$ branched 
over a Montesinos knot $k(p,q,r)$. In this paper, we express the Floer 
instanton homology groups $I_n (\Sigma (p,q,r))$, see \cite{F}, in terms of 
the Casson's $\lambda$--invariant of $\Sigma (p,q,r)$ and the signature of the 
knot $k(p,q,r)$.

\begin{theorem}\label{main}
The Floer homology group $I_n (\Sigma(p,q,r))$ is trivial 
if $n$ is odd, and is a free abelian group of the rank
$$\rank I_n(\Sigma(p,q,r)) = \frac 1 {16}\, (\,8\lambda (\Sigma(p,q,r)) -
(-1)^{n/2}\,\sign k(p,q,r))$$
if $n$ is even.
\end{theorem}

The groups $I_* (\Sigma(p,q,r))$ were studied by R.~Fintushel and R.~Stern in 
\cite{FS} (we use the Floer index convention of this paper). They 
gave an algorithm which for any given set of positive coprime integers $p,q,r$ 
computes the ranks of $I_*(\Sigma(p,q,r))$. However, their approach does not 
give a closed form formula for $I_*(\Sigma(p,q,r))$. It should be mentioned 
though that the vanishing of the groups $I_n(\Sigma(p,q,r))$ for odd $n$ 
follows from \cite{FS}; we give an alternative proof in Section 3.

In \cite{S} we introduced an invariant $\nu$, which for Brieskorn homology
spheres takes the form
\begin{equation}\label{nu}
\nu(\Sigma(p,q,r))=\frac 1 2 \sum_{n=0}^3 (-1)^{n+1} \rank I_{2n} 
(\Sigma(p,q,r)). 
\end{equation}
Note that if the ranks in (\ref{nu}) were added without plus or minus signs 
one would get the Casson invariant $\lambda(\Sigma(p,q,r))$. We conjectured in 
\cite{S} that the $\nu$--invariant equals the $\bar\mu$--invariant of 
W.~Neumann \cite{N} and L.~Siebenmann \cite{Sb} for all Seifert fibered 
homology spheres. We prove that our conjecture is true. The first step is 
given by the following theorem.

\begin{theorem}\label{conj}
For any Brieskorn homology 3--sphere $\Sigma(p,q,r)$, the $\nu$--invariant
defined by (\ref{nu}) and the $\bar\mu$--invariant coincide,
$\nu(\Sigma(p,q,r))=\bar\mu(\Sigma(p,q,r))$.
\end{theorem}

It can be easily seen that Theorem \ref{conj} is in fact equivalent to Theorem
\ref{main} due to the following two observations: first,
$\bar\mu(\Sigma(p,q,r))=1/8\cdot\sign k(p,q,r)$, see \cite{Sb}, and second, 
$I_n(\Sigma(p,q,r))=I_{n+4}(\Sigma(p,q,r))$ for all $n$. Second step in proving
our conjecture makes use of the so called splicing additivity proven for both 
the $\bar\mu$-- and $\nu$--invariants in \cite{S}. Namely,
$\bar\mu(\Sigma(a_1,\ldots,a_n))=\bar\mu(\Sigma(a_1,\ldots,a_j,p)) + 
\bar\mu(\Sigma(q,a_{j+1},\ldots,a_n))$, where $j$ is any integer between 2 and 
$n-2$,  and the integers $q=a_1\cdots a_j$ and $p=a_{j+1}\cdots a_n$ are
the products of the first $j$ and the last $(n-j)$ Seifert invariants, 
respectively. The same additivity holds for the invariant $\nu$.

\begin{theorem}\label{conj1}
For any Seifert fibered homology 3--sphere $\Sigma(a_1,\ldots,a_n)$, the
invariants $\nu$ and $\bar\mu$ coincide,
$\nu(\Sigma(a_1,\ldots,a_n))=\bar\mu(\Sigma(a_1,\ldots,a_n))$.
\end{theorem}

Another reformulation of Theorem \ref{main} establishes links between Floer 
homology and Jones polynomial. The result below follows from Theorem 
\ref{main} and D.~Mullins formula for the Casson invariant of two-fold 
branched covers \cite{Mu}.

\begin{theorem}\label{jones}
Let $\Sigma(p,q,r)$ be a Brieskorn homology sphere then 
$$\rank I_0 (\Sigma (p,q,r)) = - \frac 1 {12} \cdot \left.\frac d {dt}\right|_{t=-1} 
\ln V_k\,(t)$$
where $V_k$ is the Jones polynomial of the Montesinos knot $k = k(p,q,r)$.
\end{theorem}

Theorem \ref{conj1} has the following application to the homology 4--cobordism.
In \cite{S} and \cite{S1} we proved that $\bar\mu$ is a homology cobordism 
invariant for some classes of Seifert fibered homology 3--spheres. The
corresponding results for the $\nu$--invariant hold due to the identification
$\bar\mu=\nu$ in Theorem \ref{conj1}.

\begin{theorem} (1)\ Let $\Sigma(a_1,\ldots,a_n)$ be a Seifert fibered homology
sphere homology cobordant to zero. Then $\nu(\Sigma(a_1,\ldots,a_n))\ge 0$.

(2)\ Let $\Sigma=\Sigma(p,q,pqm\pm 1)$ be a surgery on a $(p,q)$--torus knot, 
and suppose that $\Sigma$ is homology cobordant to zero. Then $\nu(\Sigma)=0$.

(3)\ For all known Seifert fibered homology 3--spheres $\Sigma$ which are known
to be homology cobordant to zero including the lists of Casson-Harer
\cite{CH} and Stern \cite{St}, $\nu(\Sigma)=0$.
\end{theorem}
\vspace{5mm}

We will prove that $\nu(\Sigma(p,q,r))=1/8\cdot \sign k(p,q,r)$, which is
equivalent to proving Theorem \ref{main}. The idea is shortly 
as follows. In \cite{T}, C.~Taubes proved that, for any homology 3--sphere 
$\Sigma$,
\begin{equation}\label{t}
\lambda(\Sigma)= 1/2\cdot \chi (I_*(\Sigma)).
\end{equation}
The Casson's $\lambda$--invariant on the left is defined topologically using 
a Heegaard splitting of $\Sigma$ and $\SU(2)$--representation spaces. The 
number on the right is the Euler characteristic of Floer homology $I_*(\Sigma)$;
it can be defined using $\SU(2)$ gauge theory as an infinite dimensional 
generalization of the classical Euler characteristic.

Let now $\Sigma=\Sigma(p,q,r)$ be endowed with the involution $\sigma$ induced
on $\Sigma\subset {\mathbb C}^3$ by the complex conjugation. We work 
out $\sigma$--invariant versions of both invariants in (\ref{t}) in this 
particular situation. We use a $\sigma$--invariant Heegaard splitting
of $\Sigma$ and the corresponding $\sigma$--invariant $\SU(2)$--representation 
spaces to define an invariant $\lambda^{\rho}(\Sigma)$ in a manner similar to
that of A.~Casson. On the other hand, a $\sigma$--invariant gauge theory 
produces a $\sigma$--invariant Euler characteristic which we denote 
$\chi^{\rho}(\Sigma)$. Note that the latter can be defined without actually 
working out any Floer homology. Then, a $\sigma$--invariant version of the 
Taubes result (\ref{t}) is that 
\begin{equation}
\lambda^{\rho}(\Sigma(p,q,r))= 1/2\cdot \chi^{\rho}(\Sigma(p,q,r)).
\end{equation}

Our next step is to show that $\lambda^{\rho}(\Sigma(p,q,r))=1/8\cdot\sign k(p,q,r)$. 
We achieve this by pushing $\sigma$--invariant representations down to the 
quotients. We do this for all the manifolds in a $\sigma$--invariant Heegaard 
splitting and note that the push-down representations are only defined on the
complements of the fixed point sets, and they map all the meridians to
trace-free matrices in $\SU(2)$. After that we are in position to identify our
$\lambda^{\rho}$--invariant with the Casson--Lin invariant $1/8\cdot\sign
k(p,q,r)$ of \cite{Lin}. The crucial observation is that {\it all} the 
representations of $\pi_1(\Sigma(p,q,r))$ are in fact $\sigma$--invariant, and
the $\lambda^{\rho}$--invariant just counts them with signs different from
those defined by Casson for his $\lambda$--invariant.

Finally, we identify $1/2\cdot\chi^{\rho}(\Sigma)$ with $\nu(\Sigma)$ as 
follows. We choose a $\sigma$--invariant metric on $\Sigma$ and an almost 
complex structure on $\Sigma\times {\mathbb R}$ so that the involution 
$\sigma\times 1$ on $\Sigma\times {\mathbb R}$ is anti-holomorphic. At the 
level of tangent spaces, the $\pm 1$--eigenspaces of $\sigma\times 1$ are then 
isomorphic to each other and hence have the same dimension. Therefore, the 
spectral flow used in \cite{T} to define $\chi (\Sigma)$ splits in halves, and 
the $\nu$--invariant becomes an ``honest" $\sigma$--invariant Euler 
characteristic equal to $1/2\cdot \chi^{\rho}(\Sigma)$.

The paper begins with an introduction to the relevant topology of Brieskorn
homology spheres and Montesinos knots. In Section 2 the invariant
$\lambda^{\rho}(\Sigma(p,q,r))$ is introduced and the equality 
$\lambda^{\rho}(\Sigma(p,q,r))=1/8\cdot\sign k(p,q,r)$ is proven. The
invariant $\chi^{\rho}(\Sigma(p,q,r))$ is defined in Section 3; the equality
$1/2\cdot\chi^{\rho}(\Sigma(p,q,r))=\nu (\Sigma(p,q,r))$ is proven in the same 
section. Section 4 is devoted to the proof of the fact that 
$\lambda^{\rho}(\Sigma(p,q,r))=1/2\cdot\chi^{\rho}(\Sigma(p,q,r))$. Each 
section has its own introduction and is further subdivided.

It should be pointed out that most of the results in the paper hold in a more 
general situation. A description of the relevant results will appear elsewhere,
as will a more detailed version of this preprint.

\section{Topology of Brieskorn homology spheres}
Let $p, q, r$ be relatively prime positive integers greater than or equal to 2.
The Brieskorn homology sphere $\Sigma(p,q,r)$ is defined as the algebraic link
\begin{equation}\notag
\Sigma(p,q,r)=\{\, (x,y,z)\in {\mathbb C}^3\ |\ x^p+y^q+z^r=0\,\}\cap S^5_1,
\end{equation}
where $S^5_1$ is the unit sphere in ${\mathbb C}^3$. This is a smooth naturally
oriented 3--manifold with $H_*(\Sigma(p,q,r))=H_*(S^3)$. Moreover,
$\Sigma(p,q,r)$ is Seifert fibered, see \cite{NR}, with the Seifert invariants 
$\{\,b, (p,b_1), (q,b_2), (r,b_3)\,\}$ such that 
\begin{equation}\label{one}
b_1 qr+b_2 pr+b_3 pq + bpqr=1.
\end{equation}
The complex conjugation on ${\mathbb C}^3$ obviously acts on $\Sigma(p,q,r)$ as
\begin{equation}\label{sigma}
\sigma\,: \Sigma(p,q,r)\to \Sigma(p,q,r),\quad 
                    (x,y,z)\mapsto (\bar x,\bar y,\bar z).
\end{equation}
The fixed point set of this action is never empty. The quotient of
$\Sigma(p,q,r)$ by the involution $\sigma$ is $S^3$, with the branching set the
so called Montesinos knot $k(p,q,r)$, see \cite{M}, \cite{BZ}, or \cite{Sb}. 
The knot $k(p,q,r)$ can be described by the following diagram

\vspace{22mm}
\centerline{
\fbox{
\begin{picture}(250,140)
    \put(10,50)    {\fbox{\psfig{figure=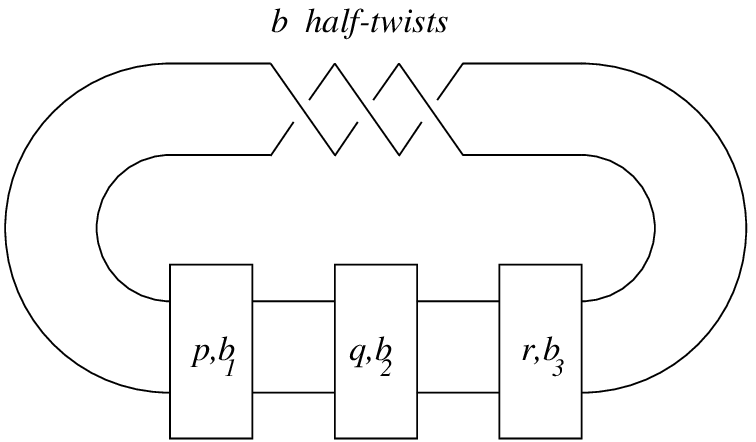}}}
    \put(125,30)    {\makebox(0,0)[t]{\bf Figure 1}}
\end{picture}
}
}

\noindent where a box with $\alpha,\beta$ in it stands for a rational
$(\alpha,\beta)$--tangle, see \cite{BZ}, Fig. 12.9. The parameters $b, (p,b_1),
(q,b_2), (r,b_3)$ are the Seifert invariants of the corresponding 
$\Sigma(p,q,r)$. According to \cite{BZ}, Theorem 12.28, these parameters
together with (\ref{one}) determine the knot $k(p,q,r)$ uniquely up to isotopy.
\vspace{3mm}

\section{The invariant $\lambda^{\rho}$}

In this section we first shortly recall the definition of Casson's
$\lambda$--invariant. Our $\lambda^{\rho}$--invariant for $\Sigma(p,q,r)$ will 
be a $\sigma$--invariant version of $\lambda$. To define it, we introduce a 
$\sigma$--invariant Heegaard splitting of $\Sigma(p,q,r)$. Then we define the
corresponding $\sigma$--invariant representation spaces and investigate closely
the representation space of $\pi_1(\Sigma(p,q,r))$. It turns out that all the 
representations in the latter space are $\sigma$--invariant. After computing
the necessary dimensions and checking the transversality condition, we define
the $\lambda^{\rho}$--invariant as an intersection number of
$\sigma$--invariant representation spaces. Finally, we check that
$\lambda^{\rho}(\Sigma(p,q,r))=1/8\cdot\sign k(p,q,r)$.
\vspace{3mm}

\noindent\textbf{1. Casson invariant.}\ Let $M$ be an oriented homology
3--sphere with a Heegard splitting $M=M_1\cup M_2$ where $M_1$ and $M_2$
are handlebodies of genus $g\ge 2$ glued along their common boundary, a
Riemann surface $M_0$. Let
$${\mathcal R}(M_i)=\Hom^*(\pi_1(M_i),\SU(2))/\ad\SU(2),\ i=0,1,2,\emptyset,$$
be the set of conjugacy classes of irreducible representations of $\pi_1 (M_i)$
is $\SU(2)$. Each ${\mathcal R}(M_i),\ i=0,1,2$, is naturally an oriented manifold. The
dimension of ${\mathcal R}(M_1)$ and ${\mathcal R}(M_2)$ is $3g-3$, and that of
${\mathcal R}(M_0)$ is $6g-6$. The inclusions $M_0\subset M_i,\ i=1,2$, induce 
embeddings ${\mathcal R}(M_i)\subset {\mathcal R}(M_0)$. The points of
intersection of ${\mathcal R}(M_1)$ with ${\mathcal R}(M_2)$ are in one-to-one
correspondence with ${\mathcal R}(M)$. If the intersection is transversal we
define the Casson's $\lambda$--invariant as
$$\lambda(M)=\frac 1 2 \sum_{\alpha\in{\mathcal R}(M)} \ep_{\alpha}$$
where $\ep_{\alpha}=\pm 1$ is a sign obtained by comparing the orientations on 
$T_{\alpha}{\mathcal R}(M_1)\oplus T_{\alpha}{\mathcal R}(M_2)$
and $T_{\alpha}{\mathcal R}(M_0)$. Note that the intersection {\it is} transversal
if $\Sigma=\Sigma(p,q,r)$. 
\vspace{3mm}

\noindent\textbf {2. The invariant Heegaard splitting.}\ We are going to 
construct an invariant Heegaard splitting of $\Sigma(p,q,r)$. Let us first fix 
notations. By $\Sigma$ we denote a Brieskorn homology sphere
$\Sigma(p,q,r)$ and by $\sigma:\Sigma\to\Sigma$ the involution constructed in
(\ref{sigma}). The projection on the quotient space will be denoted by $\pi$, so
$\pi:\Sigma\to\Sigma/\sigma=S^3$. The projection $\pi$ maps the fixed point set
$\Fix(\sigma)\subset\Sigma$ onto the Montesinos knot $k=k(p,q,r)\subset S^3$.

The knot $k\subset S^3$ can be represented as the closure of a braid $\beta$ on
$n$ strings. We think about $k$ as consisting of the braid $\beta$ and $n$
untangled arcs in $S^3$ forming its closure, see Figure 2. Let $S\subset S^3$ be 
an embedded 2--sphere in $S^3$ splitting $S^3$ in two 3--balls, $B_1$ and
$B_2$, with common boundary $S$, and such that the entire braid $\beta$ belongs 
to $\operatorname{int} B_1$. The intersection $k\cap S$ consists of $2n$ points
$P_1,\ldots,P_{2n}$. 

\vspace{30mm}
\centerline{
\fbox{
\begin{picture}(250,140)
    \put(30,50)    {\fbox{\psfig{figure=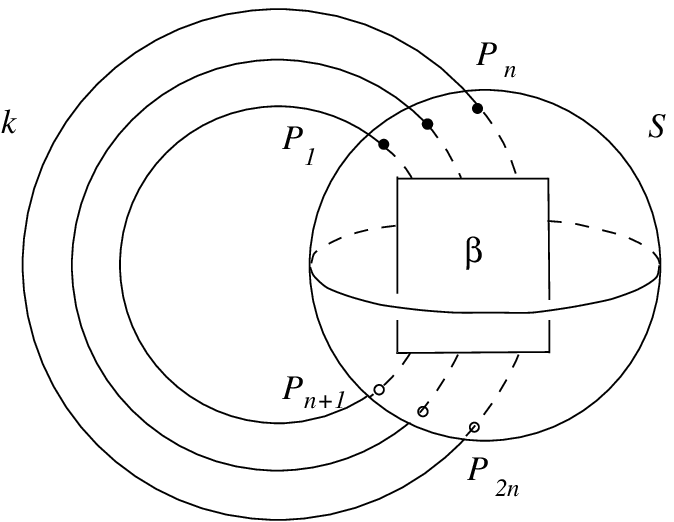}}}
    \put(125,30)    {\makebox(0,0)[t]{\bf Figure 2}}
\end{picture}
}
}

\noindent Now, we define a Heegaard splitting $\Sigma=M_1\cup_{M_0} M_2$ as 
follows:
\begin{equation}\label{heegaard}
M_1=\pi^{-1}(B_1),\quad M_2=\pi^{-1}(B_2),\quad M_0=\pi^{-1}(S).
\end{equation}
Obviously, $M_1$ and $M_2$ are handlebodies branched over the braid $\beta$ and
the $n$ untangled arcs, respectively. Their common boundary is $M_0$, which is a
closed surface branched over the points $P_1,\ldots,P_{2n}\in S$. The genus $g$
of $M_0$ can be figured out by comparing Euler characteristics,
\begin{equation}\notag
2\cdot\chi(S) = \chi(M_0) + \chi(\,\Fix(\sigma)\cap S),
\end{equation}
where $\chi(S)=2,\ \chi(\,\Fix(\sigma)\cap S)=2n$, and $\chi(M_0)=2-2g$.
Therefore, $g=n-1$. The Heegaard splitting we just defined is 
$\sigma$--invariant in the sense tnat $\sigma(M_i)=M_i$ for $i=0,1,2$.
\vspace{3mm}

\noindent\textbf {3. Representations of Brieskorn homology spheres.}\
Let $\Sigma(p,q,r)$ be a Bries\-korn homology 3--sphere, and
\begin{equation}
{\mathcal R}(\Sigma(p,q,r))=\Hom^*(\pi_1(\Sigma(p,q,r)),\SU(2))/\ad\SU(2)
\end{equation}
the space of the conjugacy classes of irreducible 
representations of its fundamental group in $\SU(2)$. In this subsection we 
are concerned with describing the space ${\mathcal R}(\Sigma (p,q,r))$ and the 
involution $\sigma^*$ induced on the representation space by $\sigma$.

We will follow \cite{FS} and first identify the group 
$\SU(2)$ with the group $S^3$ of unit quaternions in the usual way, so that
$$i = \begin{pmatrix}
      i &  0 \\
      0 & -i  
      \end{pmatrix},\quad
  j = \begin{pmatrix}
      0  & 1 \\
      -1 & 0 
      \end{pmatrix},\quad
  k = \begin{pmatrix}
      0 & i \\
      i & 0 
      \end{pmatrix}.
$$
Under this identification, the trace $\tr (A)$ of an element $A\in\SU(2)$ 
coincides with the number $2\Re A,\ A\in S^3$, and the mapping $r:S^3\to[0,\pi],\
r(A)=\arccos(\Re A)$, is a complete invariant of the conjugacy class of the 
element $A$. This conjugacy class is a copy of $S^2$ in $S^3$ unless $A = 
\pm 1$, in which case the conjugacy class consists of just one point.

The fundamental group $\pi_1(\Sigma(p,q,r))$ has the following presentation,
see \cite{BZ},
$$
\pi_1 (\Sigma (p,q,r)) = \langle x,y,z,h\ |\ h \text{ central },
  x^p = h^{-b_1}, y^q = h^{-b_2}, z^r = h^{-b_3}, xyz = h^{-b} \rangle.
$$
  
\noindent Specifying an irreducible representation $\alpha: 
\pi_1 (\Sigma (p,q,r))\to\SU(2)$ amounts to specifying a set 
$\{\alpha(h), \alpha(x), \alpha(y), \alpha(z)\}$ of unit quaternions.
In fact, we only need to specify the first three quaternions in this set
because $\alpha(z)$ will then be expressed in their terms as $\alpha(z) =
(\alpha(x)\alpha(y))^{-1}\alpha(h)^{-b}$. Since $h$ is central and the 
representation $\alpha$ is irreducible, $\alpha(h)=\pm 1$. Let us denote
$\ep_i = \alpha(h)^{b_i} = \pm 1,\ i=1,2$, and $\ep_3 = \alpha(h)^{b_3-rb}
= \pm 1$. Then the relations $x^p = h^{-b_1}, y^q = h^{-b_2}$ and $(xy)^r
= h^{b_3-rb}$ imply the following restrictions on $\alpha(x)$ and $\alpha(y)$:
$$
r(\alpha(x))=\pi\ell_1/p,\quad
r(\alpha(y))=\pi\ell_2/q,\quad
r(\alpha(x)\alpha(y))=\pi\ell_3/r,
$$
where $\ell_i$ is even if $\ep_i=1$, $\ell_i$ is odd if $\ep_i=-1$, and 
$0<\ell_1<p,\ 0<\ell_2<q,\ 0<\ell_3<r$. 

After conjugation, we may assume that $\alpha(x)=e^{\pi i\ell_1/p}$.
The quaternions $\alpha(y)$ and $\alpha(x)\alpha(y)$ should lie in their
respective conjugacy classes, $S_2=r^{-1}(\pi\ell_2/q)$ and 
$S^3=r^{-1}(\pi\ell_3/r)$. On the other hand, $\alpha(x)\alpha(y)$ lies in 
$\alpha(x)\cdot S_2$, therefore, in order for $\alpha(x)$ and $\alpha(y)$
to define a representation, the intersection $\alpha(x)\cdot S_2\cap S_3$
must be non-empty. Since $\alpha(x)\cdot S_2$ is a 2-sphere centered at
$\alpha(x)$, the intersection $\alpha(x)\cdot S_2\cap S_3$ in $S^3$ (if 
non-empty) is a circle. This circle parametrizes a whole collection of 
representations $\alpha'$ coming together with $\alpha$ such that 
$r(\alpha'(x))=r(\alpha(x)),\  r(\alpha'(y))=r(\alpha(y))$, and 
$r(\alpha'(x)\alpha'(y))=r(\alpha(x)\alpha(y))$. In fact, all these 
representations are conjugate to each other by simultaneous conjugation of 
$\alpha(x)$ and $\alpha(y)$ by the complex circle $S^1\subset S^3$. This 
can be seen from the following technical lemma.

\begin{lemma} \label{L1} Let $\alpha$ and $\beta$ be irreducible 
representations of the group $\pi_1(\Sigma(p,q,r))$ in $\SU(2)$ such 
that
\begin{enumerate}
\item [(1)] $\alpha(h)=\beta(h)$ and $\alpha(x)=\beta(x)\in{\mathbb C}$,
\item [(2)] $r(\alpha(y))=r(\beta(y))$, and
\item [(3)] $r(\alpha(x)\alpha(y))=r(\beta(x)\beta(y))$.
\end{enumerate}
Then the representations $\alpha$ and $\beta$ are conjugate to each other, 
that is there exists a unit quaternion $c$ such that
$$\beta(t)=c\cdot\alpha(t)\cdot c^{-1}\quad\text{for\ all}\quad 
                                                t\in\pi_1(\Sigma(p,q,r)).$$
Moreover, the quaternion $c$ may be chosen to be a complex number, 
$c\in {\mathbb C}$.
\end{lemma}

\begin{corollary}( Compare \cite{FS} ).
The representation space ${\mathcal R}(\Sigma(p,q,r))$ is finite.
\end{corollary}

The involution $\sigma:\Sigma(p,q,r)\to\Sigma(p,q,r)$ induces the involution 
on the fundamental group, 
\begin{gather}
\sigma_*: \pi_1 (\Sigma(p,q,r))\to\pi_1 (\Sigma(p,q,r)), \notag \\
h\mapsto h^{-1},\quad x\mapsto x^{-1},\quad y\mapsto xy^{-1}x^{-1},\quad
z\mapsto xyz^{-1}y^{-1}x^{-1},\notag
\end{gather} 
see \cite{BZ}, Proposition 12.30, which in turn induces an involution on the 
corresponding representation space ( $[\cdot,\cdot]$ stands for conjugacy 
class),
\begin{gather}
\sigma^*:{\mathcal R}(\Sigma(p,q,r))\to {\mathcal R}(\Sigma(p,q,r)),\quad
[\alpha]\mapsto [\alpha'],\label{I1} \\
\text{where}\quad \alpha'(t)=\alpha(\sigma_*(t)),\quad
t\in\pi_1(\Sigma(p,q,r)).\label{I2}
\end{gather}

\begin{lemma} \label{L2} If $\alpha': \pi_1 (\Sigma(p,q,r))\to \SU(2)$ is a
representation defined by the formula (\ref{I2}) then there exists
an element $\rho\in\SU(2)$ such that $\rho^2=-1$ and
$$\alpha'(t)=\rho\cdot\alpha(t)\cdot \rho^{-1}\quad\text{for\ all}\quad 
                                    t\in \pi_1 (\Sigma(p,q,r)).$$
The element $\rho$ is defined uniquely up to multiplication by $\pm 1$, and 
the elements $\rho$ corresponding to different representations $\alpha$ are
conjugate to each other. In particular, the action (\ref{I1}) on the space 
${\mathcal R}(\Sigma(p,q,r))$ of the conjugacy classes of irreducible 
representations of $\pi_1 (\Sigma(p,q,r))$ in $\SU(2)$ is trivial.
\end{lemma}

\begin{proof} After conjugation we may assume that $\alpha(x)$ is 
a unit complex number. Then $\alpha'(x)=\alpha(\sigma_*(x))=
\alpha(x^{-1})=\alpha(x)^{-1}$ is a complex number as well. Let us 
consider another representation, $\beta$, defined as a conjugate of 
$\alpha$ by the unit quaternion $j$,
$$\beta(t)=j^{-1}\cdot \alpha'(t)\cdot j\quad\text{for\ all}
                              \quad t\in \pi_1 (\Sigma(p,q,r)).$$
Next we want to verify that the representations $\alpha$ and $\beta$
satisfy the conditions of Lemma \ref{L1}. We have
\begin{alignat}{2}
\beta(x) &= j^{-1}\cdot \alpha'(x)\cdot j             & &                            \notag \\
         &= j^{-1}\cdot \alpha(x)^{-1}\cdot j         & &                            \notag \\
         &= j^{-1}\cdot {\overline{\alpha(x)}}\cdot j & &                            \notag \\
         &= j^{-1}\cdot j\cdot\alpha(x),              & & 
                              \quad\text{since}\quad\alpha(x)\in{\mathbb C},            \notag \\
         &= \alpha(x).                                                               \notag
\end{alignat}
Note that, for any unit quaternion $a$, we have $r(a)=r(\bar a)=r(a^{-1})$. 
Therefore,
\begin{alignat}{2}
r(\beta(y)) &=r(j^{-1}\cdot\alpha'(y)\cdot j)=r(\alpha'(y)) & & \notag \\
          &=r(\alpha(xy^{-1}x^{-1}))=r(\alpha(x)\alpha(y)^{-1}\alpha(x)^{-1})& & \notag \\
          &=r(\alpha(y)^{-1})=r(\alpha(y)), & & \notag
\end{alignat}
and similarly,
\begin{alignat}{2}
r(\beta(x)\beta(y)) &=r(j^{-1}\cdot\alpha'(x)\alpha'(y)\cdot j) 
                    =r(\alpha'(x)\alpha'(y)) & &                               \notag \\
                    &=r(\alpha(x)^{-1}\alpha(x)\alpha(y)^{-1}\alpha(x)^{-1})
                    =r((\alpha(x)\alpha(y))^{-1}) & &                          \notag \\
                    &=r(\alpha(x)\alpha(y)). & &                               \notag
\end{alignat}

Since $r(\alpha(h))$ is equal to $r(\beta(h))$ we can apply 
Lemma \ref{L1} to the representations $\alpha$ and $\beta$ 
to find a complex number $c\in {\mathbb C}$ such that $\beta(t)= c\cdot 
\alpha(t)\cdot c^{-1}$ for all $t\in \pi_1 (\Sigma(p,q,r))$. Since 
$\beta(t) = j^{-1}\cdot \alpha'(t)\cdot j = c\cdot \alpha(t)\cdot c^{-1}$,
we get that $\alpha'(t) = \rho\cdot \alpha(t)\cdot \rho^{-1}$ with $\rho=jc$.

If we apply $\sigma^*$ twice to the representation $\alpha$, we will get 
$\alpha$ again because $\sigma^*$ is an involution. Therefore, $\alpha$ is 
conjugate to itself by the element $\rho^2$. But the representation $\alpha$ 
is irreducible, therefore, $\rho^2=\pm 1$. The following easy computation shows 
that $\rho^2$ is in fact $-1$,
$$\rho^2 = jcjc = jj{\bar c}c = -|c|^2 = -1.$$

The uniqueness of $\rho$ up to $\pm 1$ follows from the irreducibility of
$\alpha$. The fact that the elements $\rho$ defined by different
representations are conjugate follows from the fact that $\tr(jc)=0$ for any 
complex number $c$.
\end{proof}
\vspace{3mm}

\noindent\textbf {4. The definition of} $\lambda^{\rho}$.\
Let $\Sigma=\Sigma(p,q,r)$ be a Brieskorn homology sphere with a
$\sigma$--invariant Heegaard splitting
\begin{equation}\label{heegaard1} 
\Sigma=M_1\cup_{M_0} M_2
\end{equation}
constructed in (\ref{heegaard}). According to Lemma \ref{L2}, for any
representation $\alpha: \pi_1(\Sigma(p,q,r))\to \SU(2)$, there exists an
element $\rho\in\SU(2)$ such that $\alpha\circ\sigma_*=\rho\alpha\rho^{-1}$
and $\rho^2=-1$. Generally speaking, $\rho$ depends on $\alpha$ but it 
follows from Lemma \ref{L2} that the elements $\rho$ corresponding to 
different $\alpha$'s are conjugate to each other. Therefore, the following
definition makes sense.

For each manifold $M_i,\ i=0,1,2$, in (\ref{heegaard1}) we define the so 
called $\rho$--invariant representation space,
\begin{multline}\label{rep}
{\mathcal R}^{\rho}(M_i)=\{\, \alpha:\pi_1 (M_i)\to\SU(2)\ |\quad 
\alpha\ \text{ is\ irreducible\ and}\\
\alpha(\sigma_*(t))=\rho\alpha(t)\rho^{-1},\quad t\in\pi_1(M_i)\, \}/\ad\SU(2)^{\rho},
\end{multline}
as the space of all irreducible $\rho$--invariant representations modulo the
adjoint action of the 1--dimensional Lie group $\SU(2)^{\rho}\subset\SU(2)$ 
consisting of all $g\in\SU(2)$ such that $\rho g = \pm g\rho$. Note that whenever 
$\rho'=\pm h\rho h^{-1},\ h\in\SU(2)$, we have that
${\mathcal R}^{\rho}(M_i)={\mathcal R}^{\rho'}(M_i)$, hence each of the spaces 
(\ref{rep}) only depends on the conjugacy class of $\rho$ and not on $\rho$ 
itself. We do not need to define a 
$\rho$--invariant representation space ${\mathcal R}^{\rho}(\Sigma(p,q,r))$
because all the representations of $\pi_1(\Sigma(p,q,r))$ in $\SU(2)$ are 
$\rho$--invariant, see Lemma \ref{L2}, and therefore 
${\mathcal R}^{\rho}(\Sigma(p,q,r))={\mathcal R}(\Sigma(p,q,r))$.

Due to van Kampen's theorem, we have the following commutative diagrams of the
fundamental groups

\[
\begin{CD}
\pi_1 (\Sigma(p,q,r))   @<<<   \pi_1 (M_1)  \\
@AAA                           @AAA         \\
\pi_1 (M_2)             @<<<   \pi_1 (M_0)
\end{CD}
\]

\vspace{3mm}

\noindent and of the $\rho$--invariant representation spaces

\[
\begin{CD}\label{cd}
{\mathcal R}(\Sigma(p,q,r))   @>>>   {\mathcal R}^{\rho} (M_1)  \\
@VVV                             @VVV                   \\
{\mathcal R}^{\rho} (M_2)     @>>>   {\mathcal R}^{\rho} (M_0)
\end{CD}
\]

\vspace{3mm}

\noindent The maps in the latter diagram are injective, so we can think about 
${\mathcal R}^{\rho}(\Sigma(p,q,r))$ as the intersection of 
${\mathcal R}^{\rho}(M_1)$ and ${\mathcal R}^{\rho} (M_2)$ inside 
${\mathcal R}^{\rho} (M_0)$. 

\begin{proposition} 
The $\rho$--invariant representation varieties ${\mathcal R}^{\rho}(M_i),\ 
i=0,1,2$, are smooth manifolds of dimensions ${\mathcal R}^{\rho}(M_i)=
2g-1,\ i=1,2$, and $\dim {\mathcal R}^{\rho}(M_0)=4g-2$ where $g$ is the genus 
of a $\sigma$--invariant Heegaard splitting (\ref{heegaard1}). The maps in the 
second diagram are embeddings. The submanifolds 
${\mathcal R}^{\rho}(M_1)$ and ${\mathcal R}^{\rho}(M_2)$ of 
${\mathcal R}^{\rho}(M_0)$ intersect transversally, and their intersection  
is in one-to-one correspondence with ${\mathcal R}(\Sigma(p,q,r))$.
\end{proposition}

The manifolds ${\mathcal R}^{\rho}(M_i),\ i=0,1,2$, can be oriented as follows.
We start with ${\mathcal R}^{\rho}(M_0)$. Since the quotient $M_0/\sigma$ is a
2--sphere and the branching set consists of $2n$ points, $n=g+1$, one can choose
generators $a_1,b_1,\ldots,a_g,b_g$ in $\pi_1(M_0)$ such that
$\sigma_*(a_i)=a_i^{-1},\ \sigma_*(b_i)=b_i^{-1}$. Then ${\mathcal R}^{\rho}
(M_0)\subset\Hom ^{\rho}(F_{2g},\SU(2))/\ad\SU(2)^{\rho}$ where $F_{2g}$ is a free
group on the generators $a_1,b_1,\ldots,a_g,b_g$, and 
\begin{multline}\notag
\Hom^{\rho}(F_{2g},\SU(2))=\{\ \alpha: F_{2g}\to\SU(2)\ |\\
  \alpha(a_i)^{-1}=\rho\alpha(a_i)\rho^{-1},\ 
  \alpha(b_i)^{-1}=\rho\alpha(b_i)\rho^{-1},\ i=1,\ldots,g\ \}.
\end{multline}

\begin{lemma}
Let $\rho\in\SU(2)$ be such that $\rho^2=-1$. Then the subset $S^{\rho}$ of
$\SU(2)$ consisting of all $a\in\SU(2)$ such that $a^{-1}=\rho a\rho^{-1}$ is
a two-dimensional sphere $S^2\subset\SU(2)$.
\end{lemma} 

\begin{proof} Any element $\rho$ of $\SU(2)$ with $\rho^2=-1$ has zero trace.
Therefore, there exists $x\in\SU(2)$ such that $\rho=xjx^{-1}$ (remember that we
identify $\SU(2)$ with the group of unit quaternions), and then the map $a\mapsto 
x^{-1}ax$ establishes a diffeomorphism between $S^j$ and $S^{\rho}$. Now, $S^j$
consists of all $b\in\SU(2)$ such that $\bar b j=jb$ and $|b|^2=1$, hence
$b=u+vi+wk$ for some real $u,v,w$ with $u^2+v^2+w^2=1$. This of course determines
a 2--sphere.  
\end{proof}

Our choice of $a_1,b_1,\ldots,a_g,b_g$ establishes a diffeomorphism
\begin{equation}\notag
\Hom^{\rho}(F_{2g},\SU(2))=(S^2)^{2g}.
\end{equation}
By choosing an orientation on $S^2$, we orient $\Hom^{\rho}(F_{2g},\SU(2))$ as a
product. This also orients ${\mathcal R}(M_0)$ by the usual ``base--fiber"
convention. Note that the orientation on ${\mathcal R}^{\rho}(M_0)$ is independent
of our choices of orientations on $S^2$ and $\SU(2)^{\rho}$ or our choice of a
basis $a_1,b_1,\ldots,a_g,b_g$.

The orientation on ${\mathcal R}^{\rho}(M_1)$ can be defined as follows. First,
we choose a basis $x_1,\ldots,x_g$ of the free group $\pi_1(M_1)$ so that
$\sigma_*(x_i)=x_i^{-1},\ i=1,\ldots,g$. Then ${\mathcal R}^{\rho}(M_1)\subset
\Hom^{\rho}(F_g,\SU(2))/\ad\SU(2)^{\rho}$, where $F_g$ is a free group on the 
generators $x_1,\ldots,x_g$, and 
\begin{multline}\notag
\Hom^{\rho}(F_g,\SU(2))=\{\ \alpha: F_{2g}\to\SU(2)\ |\
                  \alpha(\sigma_*(x))=\rho\alpha(x)\rho^{-1}\}\\
=\{\alpha: F_g\to\SU(2)\ |\ \alpha(x_i)^{-1}=\rho\alpha(x_i)\rho^{-1},\ 
                  i=1,\ldots,g\ \}.
\end{multline}
Our choice of $x_1,\ldots,x_n$ establishes a diffeomorphism 
\begin{equation}\notag
\Hom^{\rho}(F_g,\SU(2))=(S^2)^g.
\end{equation}
We choose an orientation on $S^2$ and orient $\Hom^{\rho}(F_g,\SU(2))$ as a
product. Note that this orientation changes by $(-1)^g$ when the orientation on 
$S^2$ is changed; however, it does not depend on the choice of a basis in 
$\pi_1(M_1)=F_g$. The orientation of $\Hom^{\rho}(F_g,\SU(2))$ and an orientation
of $\SU(2)^{\rho}$ orient ${\mathcal R}(M_1)$.

The $\rho$--invariant representation space ${\mathcal R}^{\rho}(M_2)$ is oriented
in a completely similar way, the orientation of $\SU(2)^{\rho}$ having already been
fixed. 

If the orientations on $S^2$ or $\SU(2)^{\rho}$ are changed, the spaces
${\mathcal R}(M_1)$ and ${\mathcal R}(M_2)$ change their orientations
simultaneously; therefore, their algebraic intersection number in 
${\mathcal R}^{\rho}(M_0)$ does not change (though {\it a priori} it may depend on
the chosen Heegaard decomposition). We define our $\lambda^{\rho}$--invariant as 
one half of this algebraic intersection number, so that 
\begin{equation}\label{lambdarho}
\lambda^{\rho} (\Sigma(p,q,r)) = 1/2\cdot\sum_{\alpha\in{\mathcal R}(\Sigma)}
\ep_{\alpha}
\end{equation}
where $\ep_{\alpha}$ equals $\pm 1$ depending on whether the orientations on the
tangent spaces $T_{\alpha}{\mathcal R}^{\rho}(M_1)\oplus 
T_{\alpha}{\mathcal R}^{\rho}(M_2)$ and $T_{\alpha}{\mathcal R}^{\rho}(M_0)$ 
agree. 

\begin{proposition}\label{well-def}
For any Brieskorn homology sphere $\Sigma(p,q,r)$ the invariant\break 
$\lambda^{\rho}(\Sigma(p,q,r))$ is well-defined, in particular, it does not depend 
on the choice of a $\sigma$--invariant Heegaard splitting. Moreover, 
$$\lambda^{\rho}(\Sigma(p,q,r)) = 1/8\cdot\sign k(p,q,r)$$
where $\sign$ stands for the knot signature.
\end{proposition}

This proposition is proven below by pushing invariant representations down 
to a knot $k(p,q,r)$ complement in $S^3$ and applying a result from \cite{Lin} 
on trace-free $\SU(2)$--representations of knot groups.

\vspace{3mm}

Thus our definition of $\lambda^{\rho}$--invariant is modelled on the Casson's 
definition of the $\lambda$--invariant; the difference is that we use 
{\it $\rho$--invariant} representation spaces. Surprisingly enough, we get an
invariant which is {\it different} from $\lambda$. For example, 
$\lambda (\Sigma(2,3,5))=\lambda (\Sigma(2,3,7))=1$ while
$\lambda^{\rho} (\Sigma(2,3,5))=-1$ and $\lambda^{\rho} (\Sigma(2,3,7))=1$.
\vspace{3mm}

\noindent\textbf {5. Casson-Lin invariant.}\ Let $B_n$ be the braid group of
rank $n$ with the standard generators $\beta_1,\ldots,\beta_{n-1}$ represented
in a free group $F_n$ on symbols $x_1,\ldots,x_n$ as follows:
\begin{alignat}{2}
\beta_i: &\quad x_i     &\quad\mapsto  &\quad x_i x_{i+1} x_i^{-1}\notag \\
         &\quad x_{i+1} &\quad\mapsto  &\quad x_i \notag \\
         &\quad x_j     &\quad\mapsto  &\quad x_j,\quad\text{if}\quad 
         j\ne i, i+1. \notag
\end{alignat}
If $\beta\in B_n$ then the automorphism of $F_n$ representing $\beta$ maps each
$x_i$ to a conjugate of some $x_j$ and preserves the product $x_1\cdots x_n$.

Let $k\subset S^3$ be a knot represented as the closure of a braid $\beta\in
B_n$. Let us fix an embedding of $k$ into $S^3$ as shown in Figure 2, the sphere
$S$ separating $S^3$ in two 3--balls, $B_1$ and $B_2$, with the braid $\beta$
inside $B_1$ and $n$ untangled arcs inside $B_2$. The fundamental group
$\pi_1(K)$ of the knot $k$ complement $K = S^3 - nbd\, (k)$, has the presentation
\begin{equation}\notag
\pi_1 (K) = \langle\, x_1,\ldots, x_n\ |\quad x_i = \beta(x_i),\ i=1,\ldots,n\,
\rangle,
\end{equation}
the generators $x_1,\ldots,x_n$ being represented by the meridians of $\beta$.

The knot complement $K$ can be now decomposed as $K=M'_1\cup_{M'_0} M'_2$
where $M'_0=S\cap K$ and $M'_i=B_i\cap K,\ i=1,2$. The manifolds $M'_1$ and
$M'_2$ are handlebodies of genus $n$, and $M'_0$ is a 2--sphere with $2n$ small
2--discs removed around the points $P_1,\ldots,P_{2n}$, see Figure 2. Due to
van Kampen's theorem, we get the following commutative diagram of fundamental
groups

\[
\begin{CD}
\pi_1 (K            )   @<<<   \pi_1 (M'_1)  \\
@AAA                           @AAA         \\
\pi_1 (M'_2)            @<<<   \pi_1 (M'_0)
\end{CD}
\]

\vspace{3mm}

\noindent where 
\begin{gather}
\pi_1 (M'_0) = \langle\, x_1,\ldots,x_n,y_1,\ldots,y_n\, |\, x_1\cdots x_n=
               y_1\cdots y_n\,\rangle, \notag \\
\pi_1 (M'_1) = \langle\, x_1,\ldots,x_n\, |\quad\rangle,\quad              
\pi_1 (M'_2) = \langle\, y_1,\ldots,y_n\, |\quad\rangle \notag             
\end{gather}
are isomorphic to free groups. To continue an analogy with the definition of the Casson
invariant, we need to define $\SU(2)$--representation spaces  of the groups 
$\pi_1(K)$ and $\pi_1(M_i),\ i=0,1,2$, and compute the corresponding
intersection number. To make things work, we impose the extra condition on the
representations that all the meridians go to trace-free matrices in $\SU(2)$.
Thus, following \cite{Lin}, we define 
\begin{alignat}{1}
\hat H_n   &=\{\,\alpha:\pi_1(M'_0)\to\SU(2)\, |\ 
            \alpha\ \text{ is\ irreducible,}\ \tr x_i=\tr y_i=0
            \,\}/\ad\SU(2), \notag \\
\hat\Gamma_{\beta} &=\{\,\alpha:\pi_1(M'_1)\to\SU(2)\, |\ 
            \alpha\ \text{ is\ irreducible,}\ \tr x_i=0
            \,\}/\ad\SU(2), \notag \\
\hat\Lambda_n &=\{\,\alpha:\pi_1(M'_2)\to\SU(2)\, |\ 
            \alpha\ \text{ is\ irreducible,}\ \tr y_i=0
            \,\}/\ad\SU(2),\notag \\
{\mathcal R}^0(K) &=\{\,\alpha:\pi_1(K)\to\SU(2)\, |\ 
            \alpha\ \text{ is\ irreducible,}\ \tr x_i=0
            \,\}/\ad\SU(2).\notag
\end{alignat}
The first three are smooth manifolds of dimensions $\dim \hat H_n=4n-6$, and
$\dim \hat\Gamma_{\beta}=\dim \hat\Lambda_n=2n-3$. The commutative diagram of
fundamental groups induces the following commutative diagram 

\[
\begin{CD}
{\mathcal R}^0 (K)   @>>>   \hat\Gamma_{\beta}     \\
@VVV                             @VVV              \\
\hat\Lambda_n        @>>>   \hat H_n
\end{CD}
\]

\vspace{3mm}

\noindent In particular, the irreducible trace-free representations of $\pi_1
(K)$ in $\SU(2)$ are in one-to-one correspondence with the intersection points
of $\hat\Lambda_n$ with $\hat\Gamma_{\beta}$ in $\hat H_n$. The manifolds 
$\hat\Lambda_n, \hat\Gamma_{\beta}$, and $\hat H_n$ are naturally oriented, see
\cite{Lin}, and (possibly after a perturbation to make the intersection
transversal) one can define
\begin{equation}\label{h}
h(k) = \sum_{\alpha\in {\mathcal R}^0 (K)} \ep'_{\alpha}
\end{equation}
where $\ep'_{\alpha}=\pm 1$ is a sign obtained by comparing the orientations on
$T_{\alpha} \hat\Lambda_n \oplus T_{\alpha} \hat\Gamma_{\beta}$ and $T_{\alpha}
\hat H_n$. The invariant $h(k)$ only depends on the knot $k$ and not on the
choices in the definition. X.-S.~Lin in \cite{Lin} further proves that $h(k) = 
1/2\cdot\sign k$.
\vspace{3mm}

\noindent\textbf {6. Representations of a knot $k(p,q,r)$ complement.}\
Throughout this subsection $\Sigma=\Sigma(p,q,r)$ is a Brieskorn homology
sphere, and $k=k(p,q,r)$ is the corresponding Montesinos knot. 

Let $E\to \Sigma$ be a (necessarily trivial) $\SU(2)$--vector bundle over
$\Sigma$. Let us fix a Riemannian metric on $\Sigma$ and consider the Banach
manifold $${\mathcal B}^* = {\mathcal A}^*/{\mathcal G}$$ of the $L^2_2$--gauge 
equivalence classes of irreducible $L^2_1$--connections in $E$. It is a 
classical result in differential geometry that the holonomy map establishes a 
one-to-one correspondence between the gauge equivalence classes of irreducible 
flat connections and the conjugacy classes of irreducible 
$\SU(2)$--representations of $\pi_1(\Sigma)$. 

The involution $\sigma$ can be lifted to a bundle endomorphism of $E$. Any 
endomorphism of $E$ clearly induces an action on ${\mathcal A}^*$ by pull-back,
and an action on ${\mathcal B}^*$ as well. Since any two liftings of $\sigma$
differ by a gauge transformation, we have a well-defined action $\sigma^*:
{\mathcal B}^*\to {\mathcal B}^*$. Denote by ${\mathcal B}^{\sigma}$ the space 
of connections invariant with respect to $\sigma^*$. 

Let $\rho\in\SU(2)$ be such that $\rho^2=-1$. The formula 
\begin{equation}\label{lift}
(x,\xi)\mapsto(\sigma(x),\rho\xi\rho^{-1})
\end{equation}
defines a lifting of $\sigma:\Sigma\to\Sigma$ on
$E$, which will again be denoted by $\rho: E\to E$. Let ${\mathcal
B}^{\rho}\subset {\mathcal B}^{\sigma}$ consist of the gauge equivalence
classes of irreducible connections $A$ in $E$ such that $\rho^* A = A$. Due to
Lemma \ref{L2}, all irreducible flat connections on $\Sigma$ belong to
${\mathcal B}^{\rho}$; in particular, ${\mathcal B}^{\rho}$ is non-empty.
The following three lemmas are easy corollaries of Propositions 1 and 17
and Theorem 18 in \cite{W}.

\begin{lemma}
For any $\rho,\rho'\in\SU(2)$ such that $\rho^2={\rho'}^2=-1,\quad 
{\mathcal B}^{\rho}={\mathcal B}^{\rho'}$. Furthermore, as a set, ${\mathcal
B}^{\rho}$ is bijective to ${\mathcal A}^{\rho}/\overline{\mathcal G}^{\rho}$
where $\overline {\mathcal G}={\mathcal G}/\pm 1$ and ${\mathcal A}^{\rho}=
\{\,A\in{\mathcal A}^*\ |\ \rho^* A=A\,\}$.
\end{lemma}

The quotient of $\Sigma(p,q,r)$ by $\sigma$ is $S^3$. Since $\rho\ne \pm 1$, it
is impossible to define the quotient bundle of $E$ over $S^3$. One can though
define it away from the image of the fixed point set of $\sigma$, that is, on
the knot complement $K=S^3\setminus k$. Given an irreducible flat
$\rho$--invariant $\SU(2)$--connection $A\in {\mathcal A}^{\rho}$, its
push-down $A'$ is an irreducible flat $\SO(3)$--connection on $K$. In other
words, $A'$ is an irreducible flat $\SO(3)$--connection singular along $k$ in
the sense of P.~Kronheimer and T.~Mrowka, see \cite{KM}.

\begin{lemma} The $\SO(3)$--connection $A'$ has holonomy 1/2 around $k$.
\end{lemma}

According to $H^1(K;{\mathbb Z}/2)={\mathbb Z}/2$ and $H^2(K;{\mathbb Z}/2)=0$,
there are two different ways to lift the $\SO(3)$--connection $A'$ to a
(singular) $\SU(2)$--connection. Both $\SU(2)$--liftings have holonomy 1/4
around $k$; in other words, their holonomy around $k$ is trace--free.

\begin{lemma}\label{1-2} 
Through the push-down map, the representation space
${\mathcal R}(\Sigma)$ is in one-to-two correspondence with the representation
space ${\mathcal R}^0 (K)$ of irreducible trace--free $\SU(2)$--representations
of $\pi_1(K)$.
\end{lemma}

Similar one-to-two identifications through push-down hold for the pairs 
${\mathcal R}^{\rho}(M_0)$ and $\hat H_n$, $\ {\mathcal R}^{\rho}(M_1)$ and 
$\hat\Gamma_{\beta}$, $\ {\mathcal R}^{\rho} (M_2)$ and $\hat\Lambda_n$, where
$g=n-1$ is the genus of the $\sigma$--invariant Heegaard splitting, see 
Figure 2.

\vspace{2mm}

\begin{proof}[Proof of Proposition \ref{well-def}.] Remember that throughout
this subsection, $\Sigma=\Sigma(p,q,r)$ and $k=k(p,q,r)$. The invariants
$\lambda^{\rho}(\Sigma)$ and $h(k)=1/2\cdot \sign k$ were defined in
(\ref{lambdarho}) and (\ref{h}) as an algebraic count of points in the
corresponding (finite) representation spaces, ${\mathcal R}(\Sigma)$ and
${\mathcal R}^0(K)$. These spaces are in one-to-two correspondence by Lemma
\ref{1-2}. A routine comparison of the orientations shows that $\ep_{\alpha}
=\ep'_{\alpha'}$ where $\alpha'$ is a push-down of $\alpha$. Therefore,
\begin{alignat}{1}
\lambda^{\rho}(\Sigma) &= 1/2\cdot \sum \ep_{\alpha} \notag \\
                       &= 1/4\cdot \sum \ep'_{\alpha'} \notag \\
                       &= 1/8\cdot \sign k. \notag
\end{alignat}
This completes the proof of Proposition \ref{well-def}.
\end{proof}

\vspace{3mm}

\section{Gauge theory for Brieskorn homology spheres}

Our next step is to give an ``analytic" definition of the invariant 
$\lambda^{\rho}$. In \cite{T}, C.~Taubes gave a gauge-theoretical
interpretation of the Casson's $\lambda$--invariant in terms of (what later
became known as) Floer homology. Our gauge-theoretical interpretation of 
$\lambda^{\rho}$ will follow the same lines with the only difference that 
everything will be ``$\rho$--invariant".
\vspace{3mm}

\noindent\textbf {1. Casson invariant via gauge theory.}\ Let us recall shortly 
the Taubes construction \cite{T} for a homology 3--sphere $\Sigma$.
Let $E\to\Sigma$ be a trivial $\SU(2)$--bundle, and ${\mathcal B}^*=
{\mathcal A}^*/{\mathcal G}$ the space of gauge equaivalence classes of irreducible
connections in $E$. It is a classical differential geometry result that flat
connections in $E$ (modulo gauge equivalence) are in one-to-one correspondence with
$\SU(2)$--representations of $\pi_1(\Sigma)$ (modulo conjugation). Given two
irreducible representations $\alpha,\beta$, we may think of them as points in 
${\mathcal A}^*$ or ${\mathcal B}^*$.

For any $A\in {\mathcal A}^*$, define the following elliptic differential operator
\begin{equation}\label{KA}
K_A=
\begin{vmatrix}
0     &   d^*_A  \\ 
d_A   &    *d_A  
\end{vmatrix}
:(\Omega^0\oplus\Omega^1)(\Sigma,\su(2))\to(\Omega^0\oplus\Omega^1)(\Sigma,\su(2)),
\end{equation}
where $d_A$ stands for the operator of covariant differentiation with respect to
$A$, and $*$ is the Hodge operator associated with a Riemannian metric on $\Sigma$.
In a proper Sobolev completion of $(\Omega^0\oplus\Omega^1)(\Sigma,\su(2))$, the
operator $K_A$ is Fredholm. Let $A(t),\ 0\le t\le 1$, 
be a path in ${\mathcal A}^*$ 
connecting $\alpha$ with $\beta$, so that $A(0)=\alpha, A(1)=\beta$. Associated
with $A(t)$ is a path of Fredholm operators, $K_{A(t)}$. Let us assume that 
$\ker K_{A(0)} = \ker K_{A(1)}$, as is the case for a Brieskorn homology sphere
$\Sigma=\Sigma(p,q,r)$. We define $sf(\alpha,\beta)$ as the spectral flow of the 
path $K_{A(t)}$ between $\alpha$ and $\beta$, see \cite{APS}. This number is well-defined modulo 8, 
and Taubes defines an infinite dimensional generalization $\chi(\Sigma)$ of the
Euler characteristic as
\begin{equation}\label{taubes}
\chi (\Sigma) = \ep_{\alpha}\cdot\sum_{\beta\in{\mathcal R}(\Sigma)} 
                                              (-1)^{sf(\alpha,\beta)}
\end{equation}
where $\ep_{\alpha}$ is figured out from the spectral flow between the trivial
connection $\theta$ and $\alpha$ (the number $\chi (\Sigma)$ turns out to be
independent of $\alpha$). It is proven in \cite{T} that $\lambda(\Sigma)=
1/2\cdot \chi(\Sigma)$. Later, A.~Floer showed in \cite{F} that there are 
well--defined ${\mathbb Z}/8$--graded instanton homology groups $I_*(\Sigma)$
such that $\chi(\Sigma)=\chi(I_*(\Sigma))$.

Let now $\Sigma=\Sigma(p,q,r)$ be a Brieskorn homology sphere. R.~Fintushel and
R.~Stern showed in \cite{FS} (see also Proposition \ref{even-ind} below) that 
for any pair $\alpha,\beta\in{\mathcal R}(\Sigma(p,q,r))$, 
\begin{equation}\notag
sf(\alpha,\beta)=0\mod 2,
\end{equation}
and that all the signs involved in computing $\chi(\Sigma(p,q,r))$ are
positive, so $\chi$ simply counts the irreducible representations. 
\vspace{3mm}

\noindent\textbf {2. Floer index for Brieskorn homology spheres.}\  
Let $\alpha$ and $\beta$ be irreducible flat connections in a trivial
$\SU(2)$--bundle $E$ over $\Sigma=\Sigma(p,q,r)$. By pull-back, we can extend 
$E$ to a trivial bundle (which we also denote by $E$) over the infinite 
cylinder $\Sigma\times {\mathbb R}$. Let us choose a $\sigma$--invariant
Riemannian metric on $\Sigma$, and the corresponding product metric on
$\Sigma\times {\mathbb R}$. Let further $\rho: E\to E$ be the lifting of 
$\sigma$ defined in (\ref{lift}), and $A(t)\in {\mathcal A}^{\rho}$ a path of
connections forming an invariant connection $A$ over $\Sigma\times {\mathbb R}$ 
vanishing in the ${\mathbb R}$--direction and equal 
to respectively $\alpha$ and $\beta$ near the ends of $\Sigma\times 
{\mathbb R}$. Denote by $d_A$ the operator of covariant differentiation with 
respect to $A$.

The (relative) Floer index of the pair $(\alpha,\beta)$ equals the Fredholm 
index of the following elliptic complex, see \cite{APS}, 
\begin{equation}\label{E1}
\Omega^0 (\Sigma\times {\mathbb R}, \ad E) \xrightarrow{d_A} 
\Omega^1 (\Sigma\times {\mathbb R}, \ad E) \xrightarrow{d^-_A} 
\Omega^2_+ (\Sigma\times {\mathbb R}, \ad E) 
\end{equation}
where $\ad E$ is the adjoint bundle of $E$ over $\Sigma\times {\mathbb R}$, and
$d^-_A = P_- \circ d_A$ where $P_-$ is the projection onto the self dual
forms with respect to the fixed product metric on $\Sigma\times {\mathbb R}$. The
Sobolev norms on the spaces in (\ref{E1}) are specified as usual, see \cite{F}.

The product metric and the connection $A$ on $\Sigma\times {\mathbb R}$ are
$\rho$--invariant, therefore, one can define the following $\rho$--invariant 
elliptic subcomplex of (\ref{E1}),
\begin{equation}\label{E2}
\Omega^0 (\Sigma\times {\mathbb R}, \ad E)^{\rho} \xrightarrow{d_A} 
\Omega^1 (\Sigma\times {\mathbb R}, \ad E)^{\rho} \xrightarrow{d^-_A} 
\Omega^2_+ (\Sigma\times {\mathbb R}, \ad E)^{\rho} 
\end{equation}

\begin{proposition}\label{even-ind} 
For any irreducible flat connections $\alpha,\beta$ on $\Sigma=\Sigma(p,q,r)$, 
the index of the elliptic complex (\ref{E1}) is even. In fact, it equals twice 
the index of the $\rho$--invariant complex (\ref{E2}).
\end{proposition}

\noindent In order to prove the proposition we need the following two technical lemmas.

\begin{lemma}\label{L10} There exist a Riemannian metric on $\Sigma(p,q,r)$ 
and an almost complex structure $J$ on $\Sigma(p,q,r)\times {\mathbb R}$ such that 
\begin{enumerate}
\item [(1)] The product metric on $\Sigma\times {\mathbb R}$ is $(\sigma\times
      1)$--invariant;
\item [(2)] $J$ is compatible with this product metric; and
\item [(3)] The involution $\sigma\times 1$ is anti-holomorphic with respect to $J$,
      that is $(\sigma\times 1)_* J = - J(\sigma\times 1)_*$ on the tangent
      bundle.
\end{enumerate}
\end{lemma}

\begin{proof} Let us consider the algebraic variety
$$V(p,q,r)=\{\,(x,y,z)\in {\mathbb C}^3\ |\ x^p+y^q+z^r=0\,\}\subset {\mathbb C}^3,$$
which is a non-singular complex surface except perhaps at the origin. The 
homology sphere $\Sigma(p,q,r)$ in question is the intersection of $V(p,q,r)$ 
with the unit 5-sphere $S^5_1$,
$$\Sigma(p,q,r)=V(p,q,r)\cap S^5_1,$$
and the variety $V(p,q,r)$ is in fact a cone over $\Sigma(p,q,r)$ with the 
vertex at the origin. 

Let $V^0(p,q,r)$ be the variety $V(p,q,r)$ with the origin removed. The Riemann metric 
induced on $V^0(p,q,r)$ from the standard flat metric on ${\mathbb C}^3$ is a cone 
metric given by 
\begin{equation}\label{metric1}
ds^2=dr^2+r^2d\theta^2
\end{equation}
in the spherical coordinates $(r,\theta)$. Here $d\theta^2=\sum h_{ij} 
d\theta_i d\theta_j$ is a metric on $\Sigma(p,q,r)$, and $r$ is the distance 
from the origin (in ${\mathbb C}^3$). Obviously, both the complex conjugation 
$\sigma'$ and the almost complex structure $J$ induced on $V^0(p,q,r)$ from 
${\mathbb C}^3$ are compatible with the metric (\ref{metric1}).

Let us now form the conformally equivalent metric
\begin{equation}\label{metric2}
ds^2/r^2 = dr^2/r^2 + d\theta^2
\end{equation}
on $V^0(p,q,r)$. The substitution $r=e^{-\tau}$ gives coordinates in which
$ds^2/r^2$ is the standard product metric
\begin{equation}\label{metric3}
d\tilde s^2 = d\tau^2 + d\theta^2
\end{equation}
on the cylinder $\Sigma(p,q,r)\times {\mathbb R}$. Since the metrics
(\ref{metric1}) and (\ref{metric3}) are conformally equivalent, both $\sigma'$
and $J$ are compatible with the product metric (\ref{metric3}). Moreover,
$\sigma'$ preserves all the spheres $S^5_r,\ 0 < r < \infty$, in ${\mathbb C}^3$,
therefore, it has the form $\sigma'=\sigma\times 1$. Finally, the
involution $\sigma\times 1$ is antiholomorphic with respect to $J$ because this
is the case in ${\mathbb C}^3$.
\end{proof}

\begin{lemma}\label{L11} There exists a differential operator $\tilde d^-_A:
\Omega^1 (\Sigma\times {\mathbb R}, \ad E)\to \Omega^2_+ (\Sigma\times {\mathbb R}, 
\ad E)$ such that 
\begin{enumerate}
\item [(1)] The difference $\tilde d^-_A - d^-_A$ is a compact operator,
\item [(2)] $\coker \tilde d^-_A = 0$, and
\item [(3)] The operator $\tilde d^-_A$ is $\rho$--invariant.
\end{enumerate}
\end{lemma}

\begin{proof}
This is essentially the transversality result from \cite{F}, Proposition 2c.2. 
The key difference is that we have to bound ourselves to invariant 
perturbations which, generally speaking, may be not generic among the 
perturbations used by Floer. Following \cite{F}, Section 1b, consider a
collection $S(m)$ of $m$ circles smoothly embedded in ${\mathbb R}^3$ which
intersect precisely at the origin and have the same tangent there. Let
\begin{equation}\notag
\gamma: \bigvee^m_{i=1} S^1_i\times D^2 \to \Sigma,
\end{equation}
be a map which restricts to smooth embeddings $\gamma_{\theta}: S(m)\to \Sigma$
for each $\theta\in D^2$ and $\gamma_i:S^1_i\times D^2\to \Sigma$ for each $i$.
By adding extra circles, one can choose $\gamma$ so that its image in $\Sigma$
will be $\sigma$--invariant. Now, $\gamma$ defines a family of holonomy maps
\begin{gather}\notag
\gamma_{\theta}: {\mathcal B}\to L_m = \SU(2)^m/\ad\SU(2),\\
A\mapsto (\hol_A (\gamma_1),\ldots,\hol_A (\gamma_m))
\end{gather}
Let $C^{\infty}_{sym}(L_m,{\mathbb R})$ be defined as the set of smooth
$\ad\SU(2)$--invariant real valued functions on $\SU(2)^m$ with the
additional property that for any $h\in C^{\infty}_{sym}(L_m,{\mathbb R})$,
\begin{equation}\notag
h(x_{\tau(1)},\ldots,x_{\tau(m)})=h(x_1,\ldots,x_m)
\end{equation}
for any permutation $\tau$ on $m$ symbols. For example, one can start with an
arbitrary $\ad\SU(2)$--invariant function on $\SU(2)^m$ and take $h$ to be its
symmetrization.

Given such $\gamma$ and $h$, we define the function
\begin{equation}\label{hg}
h_{\gamma}:{\mathcal B}\to {\mathbb R},\quad
h_{\gamma}(A) = \int_{D^2} h(\gamma_{\theta}(A))\,d^2\theta
\end{equation}
where $d^2\theta$ is a smooth compactly supported volume form on the interior
of $D^2$. The following argument shows that this function is invariant with 
respect to the induced $\sigma^*$--action on ${\mathcal B}$:
\begin{multline}\notag
h_{\gamma} (\sigma^* A) = \int_{D^2} h(\,\hol_{\sigma^* A} (\gamma_1),\ldots,\hol_{\sigma^* A} (\gamma_m))\,d^2\theta \\
                        = \int_{D^2} h(\,\hol_A (\sigma(\gamma_1)),\ldots,\hol_A (\sigma(\gamma_m)))\,d^2\theta\qquad \\
                        \qquad = \int_{D^2} h(\,\hol_A (\gamma_{\tau(1)}),\ldots,\hol_A(\gamma_{\tau(m)}))\,d^2\theta \\
                        = \int_{D^2} h(\,\hol_A (\gamma_1),\ldots,\hol_A (\gamma_m))\,d^2\theta  = h_{\gamma} (A).    
\end{multline}
The class of maps we defined in (\ref{hg}) is large enough for the set 
$$\{\,\grad h_{\gamma}(A)\ |\ h_{\gamma}\ \text{ defined\ by\ (\ref{hg})}\,\}\subset T_A {\mathcal B}^*$$
to be dense in the tangent space $T_A {\mathcal B}^*$. Now the proof of Proposition 
2c.2 in \cite{F} goes through, which provides a function $h_{\gamma}$ such that a
linearization of its gradient is a compact invariant perturbation of the operator 
$d^-_A$ with the desired properties.
\end{proof}

\begin{proof}[Sketch of the proof of Proposition \ref{even-ind}.] 
In Lemma \ref{L10}
we defined a $\sigma$--invariant Riemannian metric on $\Sigma(p,q,r)$ and an 
almost complex structure $J$ on the manifold $\Sigma(p,q,r)\times {\mathbb R}$
compatible with the product metric such that the involution $\sigma\times 1:
\Sigma\times {\mathbb R}\to \Sigma\times {\mathbb R}$ is anti-holomorphic with 
respect to $J$. We can use $J$ to identify the $\pm 1$--eigenspaces of the 
involution induced by $\sigma\times 1$ on the space of 1--forms in (\ref{E1}), 
and make sure that (possibly after perturbation provided by Lemma \ref{L11}) 
$\ker d_A = 0$ and $\coker d^-_A = 0$. This splits the index of the complex 
(\ref{E1}) and hence the corresponding spectral flow in halves.
\end{proof}

\begin{proposition} ( Compare \cite{FS} ). Let $\Sigma(p,q,r)$ be a Brieskorn
homology sphere. Then the Floer homology groups $I_n (\Sigma(p,q,r))$ vanish 
for $n$ odd. 
\end{proposition}

\begin{proof} We only need to fix the sign $\ep_{\alpha}$ in (\ref{taubes}). It
can be achieved by the requirement that for $\ep_{\alpha}=(-1)^{\eta(\alpha)}$,
\begin{equation}\notag
-3-\eta (\alpha) = \ind D^+_A - 3\,(1+b^+_2)(W),
\end{equation}
where $W$ is any smooth simply connected non-compact 4--manifold with a single 
end of the form $\Sigma(p,q,r)\times {\mathbb R}_+$, $\ A$ is any 
$\SU(2)$--connection on $W$ with limiting value $\alpha$, and
\begin{equation}\notag
D^+_A = d_A\oplus d^+_A\,: \Omega^1(W,\su(2))\to \Omega^0(W,\su(2))\oplus
\Omega^2_+(W,\su(2))
\end{equation}
is the ASD--operator, compare \cite{DFK}. Let $W$ be the Milnor fiber of the
singularity of $f^{-1}(0)$ where $f(x,y,z)=x^p+y^q+z^r$ so that 
$\p W=\Sigma(p,q,r)$, see \cite{Mi}. The involution (\ref{sigma}) extends to an
antiholomorphic involution on $W$, and an argument similar to that in the proof
of Proposition \ref{even-ind} proves that $\ind D^+_A$ is even. On the other
hand, one can easily show that modulo 2, $b^+_2 = 1/2\cdot b_2(W) = 1/2\cdot
(p-1)(q-1)(r-1) = 0$. Therefore, $\ep_{\alpha}=1$, and the Floer
homology groups of $\Sigma(p,q,r)$ vanish in odd dimensions.
\end{proof}
\vspace{3mm}

\noindent\textbf {3. The definition of $\chi^{\rho}$.}\ Let $\rho : E\to E$ be
the lifting of $\sigma :\Sigma(p,q,r)\to \Sigma(p,q,r)$ defined by (\ref{lift}).
It induces an involution $\rho^*$ on the $\su(2)$--valued differential forms on
$\Sigma(p,q,r)$. Therefore, for any $A\in {\mathcal A}^{\rho}$, the
differential operator $K_A$, see (\ref{KA}), can be restricted to the
$+1$--eigenspaces of the involution $\rho^*$ to give a new elliptic
differential operator $K^{\rho}_A$,
\begin{equation}\notag
K^{\rho}_A=
\begin{vmatrix}
0     &   d^*_A  \\ 
d_A   &    *d_A  
\end{vmatrix}
:(\Omega^0\oplus\Omega^1)^{\rho}(\Sigma,\su(2))\to
 (\Omega^0\oplus\Omega^1)^{\rho}(\Sigma,\su(2)).
\end{equation}
For any pair $\alpha,\beta\in{\mathcal R}(\Sigma(p,q,r))$, let us denote by 
$sf^{\rho}(\alpha,\beta)$ the spectral flow of a family of operators
$K^{\rho}_{A(t)}$, where $A(t)$ is a path in ${\mathcal A}^{\rho}$ connecting
$\alpha$ to $\beta$. The number $sf^{\rho}(\alpha,\beta)$ is well-defined
modulo 4. The invariant $\chi^{\rho}(\Sigma(p,q,r))$ is now defined by the
formula
\begin{equation}\label{taubes1}
\chi^{\rho} (\Sigma) = \ep_{\alpha}\cdot\sum_{\beta\in{\mathcal R}(\Sigma)} 
                                              (-1)^{sf^{\rho}(\alpha,\beta)},
\end{equation}
compare to (\ref{taubes}). The following result is an easy consequence of
Proposition \ref{even-ind}.
\begin{proposition}\notag
Let $\Sigma=\Sigma(p,q,r)$ be a Brieskorn homology sphere. Then $1/2\cdot 
\chi^{\rho}(\Sigma) = \nu (\Sigma)$ where the invariants $\chi^{\rho}(\Sigma)$
and $\nu (\Sigma)$ are defined by (\ref{taubes1}) and (\ref{nu}).
\end{proposition}
\vspace{3mm}

\section{The invariant $\lambda^{\rho}$ via gauge theory}
\vspace{2mm}

\begin{proposition}
Let $\Sigma=\Sigma(p,q,r)$ be a Brieskorn homology sphere, and $\lambda^{\rho}$
and $\chi^{\rho}$ the invariant defined by (\ref{lambdarho}) and
(\ref{taubes1}). Then $\lambda^{\rho}(\Sigma) = 1/2\cdot \chi^{\rho}(\Sigma)$.
\end{proposition}

This result is a straightforward application of the Taubes
argument in \cite{T} in our $\rho$--invariant setting.

\end{document}